\documentclass[12pt,psamsfonts]{amsart}

\newtheorem{anyprop}{Anyprop}[section]

\newtheorem{theorem}[anyprop]{Theorem}
\newtheorem{lemma}[anyprop]{Lemma}
\newtheorem{proposition}[anyprop]{Proposition}

\theoremstyle{definition}

\newtheorem{definition}[anyprop]{Definition}

\newtheorem{remark}[anyprop]{Remark}

\newtheorem{theoremintro}{Theorem}

\newcommand{\komdots}{ , \ldots , }

\newcommand{\subsetdots}{ \subset \ldots \subset }

\newcommand{\NN}{\mathbb{N}}
\newcommand{\ZZ}{\mathbb{Z}}

\newcommand  {\shF}     {\mathcal{F}}

\newcommand  {\shL}     {\mathcal{L}}

\newcommand  {\shS}     {\mathcal{S}}
\newcommand  {\shT}     {\mathcal{T}}

\newcommand  {\shQ}     {\mathcal{Q}}

\newcommand  {\fom}     {\mathfrak{m}}

\newcommand  {\dual}    {\vee}

\newcommand  {\im}      {\operatorname{im}}

\newcommand  {\lra}     {\longrightarrow}

\newcommand  {\modu}     {\operatorname{mod}}

\renewcommand{\O}       {\mathcal{O}}

\newcommand  {\Proj}    {\operatorname{Proj}}

\newcommand  {\ra}      {\rightarrow}

\newcommand  {\rk}    {\operatorname{rk}}

\newcommand  {\Syz}     {\operatorname{Syz}}

\newcommand{\numiii}{\renewcommand{\labelenumi}{(\roman{enumi})}}

\newcommand{\sumceilm}{\sum_{m= \lceil q \nuu \rceil}^{\lceil q \nuo \rceil -1}}

\newcommand{\sumceil}[3]{ \sum_{m= \lceil q #1 \rceil}^{\lceil q #2 \rceil -1} h^0( #3^q(m))}

\newcommand{\sumceilnu}[1]{ \sum_{m= \lceil q \nuu \rceil}^{\lceil q \nuo \rceil -1}\! h^0( #1^q(m))}

\newcommand{\stacklra}[1]{ \stackrel{ #1 }{\lra} }

\newcommand{\nuo}{\rho}
\newcommand{\nuu}{\sigma}
\newcommand{\nus}{\nu}
\newcommand{\nut}{\nu_t}
\newcommand{\nuk}{\nu_k}

\newcommand{\pero}{\pi}
\newcommand{\peru}{\delta}

\newcommand{\pers}{\epsilon}
\newcommand{\pert}{\pi}

\newcommand{\perk}{\pi_k}

\newcommand{\hkf}{\varphi }
\newcommand{\length} {\lambda}
\newcommand{\mumu}{\mu_{HK}}

\theoremstyle{remark}

\numberwithin{equation}{section}

\usepackage{amscd}
\usepackage{amssymb}

\def\mydate{\number\day\space\ifcase\month \or January\or February\or March\or April\or May\or
June\or July\or August\or September\or October\or November\or
December\fi \space\number\year}

\begin{document}

\title[Hilbert-Kunz function]
{The Hilbert-Kunz function in graded dimension two}

\author[Holger Brenner]{Holger Brenner}
\address{Mathematische Fakult\"at, Ruhr-Universit\"at Bochum, 
               44780 Bochum, Germany}
\email{Holger.Brenner@ruhr-uni-bochum.de}





\renewcommand{\arraystretch}{10}

\begin{abstract}
Let $R$ denote a two-dimensional normal standard-graded $K$-domain
over the algebraic closure $K$ of a finite field of characteristic $p$, and let $I\subset R$
denote a homogeneous $R_+$-primary ideal.
We prove that the Hilbert-Kunz function of $I$ has the form
$\hkf (q) = e_{HK}(I)  q^{2} + \gamma (q)$
with rational Hilbert-Kunz multiplicity $e_{HK}(I)$
and an eventually periodic function $\gamma (q)$.
\end{abstract}

\maketitle


\section*{Introduction}

Suppose that $(R, \fom)$ is a local Noetherian or a standard-graded ring of dimension $d$ containing a field
$K$ of positive characteristic $p$. Let $I$ denote an
$\fom$-primary ideal and set $I^{[q]}=(f^{q}: f \in I)$, $q=p^{e}$.
The function $e \mapsto \length ( R/I^{[p^{e}]}) $, where $\length $ denotes the length,
is called the Hilbert-Kunz function
of the ideal $I$ and was first considered by Kunz in \cite{kunzpositive}.
Monsky showed in \cite{monskyhilbertkunz} that this function
has the form (we write $q$ for the argument, not $e$)
$$  \hkf(q)=e_{HK}(I) q^{d} +O(q^{d-1})  \, ,$$
where $e_{HK}(I)$ is a positive real number
called the Hilbert-Kunz multiplicity of the ideal.
It is conjectured that the Hilbert-Kunz multiplicity is always a rational number.

In \cite{hunekemcdermottmonsky}, Huneke, McDermott and Monsky studied further
the Hilbert-Kunz function showing that
$$ \hkf(q) = e_{HK}(I) q^{d} + \beta q^{d-1} +O(q^{d-2})  \, $$
holds with another real number $\beta$ under the condition that $R$ is normal and
excellent with a perfect residue field
(\cite[Theorem 1.12]{hunekemcdermottmonsky}).

In this paper we want to investigate the Hilbert-Kunz-function and in particular the $O(q^{d-2})$-term
in the case of a two-dimensional normal standard-graded $K$-domain over an algebraically closed field $K$
of positive characteristic $p$.
In this case we have recently shown that the Hilbert-Kunz multiplicity is a rational number
(see \cite[Theorem 3.6]{brennerhilbertkunz}, the rationality for $I=\fom$  was obtained
independently by Trivedi in \cite{trivedihilbertkunz}), and Monsky announced
in \cite[Remark 2.5]{hunekemcdermottmonsky} a proof that the second coefficient $ \beta =0$
vanishes. Our main result is the following theorem.

\begin{theoremintro}
\label{hilbertkunzfunctionperiodicintro}
Let $K$ denote the algebraic closure of a finite field of characteristic $p$.
Let $R$ denote a normal two-dimensional standard-graded
$K$-domain and let $I$ denote a homogeneous
$R_+$-primary ideal.
Then the Hilbert-Kunz function of $I$ has the form
$$ \hkf (q) =e_{HK}(I)\, q^{2} + \gamma (q) \, ,$$
where the Hilbert-Kunz multiplicity $e_{HK}(I)$ is a rational number and where
$\gamma(q)$ is an eventually
periodic function.
\end{theoremintro}

This result has been expected by several people,
but proves are known so far only in the case of the maximal ideal in a regular ring (\cite{contessa})
or for certain cones over elliptic curves (\cite{buchweitzchenhilbertkunz}, \cite{monskyellipticchartwo},
\cite{fakhruddintrivedi}).
In these known cases the condition that $K$ is the algebraic closure of a finite field is not needed.

We give an overview for the argument for this result and of this paper as a whole.
Write $I=(f_1 \komdots f_n)$ with homogeneous ideal generators
$f_i$ of degree $d_i = \deg (f_i)$. We shall use the short exact sequence
$$ 0 \lra \Syz (f_1^q \komdots f_n^q)(m) \lra \bigoplus_{i=1}^n \O(m-qd_i)
\stackrel{f_1^q \komdots f_n^q}{\lra} \O(m) \lra 0 \, $$
on the smooth projective curve $Y= \Proj R$
to compute
$$ \length  ((R/I^{[q]})_m)
= h^0 (\O(m)) - \sum_{i=1}^n  h^0(\O(m-qd_i ))+ h^0(\Syz (f_1^q \komdots f_n^q)(m) ) \, . $$
This gives the Hilbert-Kunz function by summing over $m$,
which is a finite sum.
The syzygy bundle $\Syz (f_1^q \komdots f_n^q)(m)$ is a locally free sheaf on $Y$ and
we have to compute its global sections for varying $q$ and for $m$ running in certain ranges.

It is natural and helpful to consider more generally an arbitrary locally free sheaf $\shS$ on a smooth projective
curve $Y$ over an algebraically closed field $K$ of positive characteristic endowed
with a fixed very ample invertible sheaf
$\O_Y(1)$.
Then we have to understand the global sections
$ H^0(Y, \shS^q(m))$, where $\shS^q$ denotes the pull-back under the $e-$th absolute
Frobenius morphism, $q= p^{e}$.
The appropriate object to study here is the expression
$$ \sumceilnu \shS $$
for rational numbers $\nuu$ and $\nuo$. The reason for this setting is that is allows us to do induction
on the strong Harder-Narasimhan filtration of $\shS$. This is the filtration
$ \shS_1^q \subsetdots \shS_t^q = \shS^q$ such that the quotients $\shS_k^q/ \shS_{k-1}^q$
are strongly semistable of decreasing slopes.
Strongly semistable means that every Frobenius pull-back is again semistable.
Such a filtration exists and is stable for $q \gg 0$ due to a Theorem of Langer.
Using this we can reduce many questions to the case where $\shS$ itself is strongly semistable.

It turns out that the above expression is related to what we call the Hilbert-Kunz slope of $\shS$
(section \ref{hilbertkunzslopesection}).
This is by definition $\mumu(\shS) = \sum_{k=1}^t r_k \bar{\mu}_k^2$,
where the rational numbers
$\bar{\mu}_k = \mu  ( \shS^q_k /\shS_{k-1}^q)/q$
come from the strong Harder-Narasimhan filtration.
We get after some preparatory work in section \ref{dimensionsectionssection} the following formula (Theorem \ref{sectionformulatheorem}).

\begin{theoremintro}
\label{sectionformulatheoremintro}
Let $\shS$ denote a locally free sheaf on a smooth projective curve $Y$ over an algebraically closed field $K$ of positive characteristic $p$.
Let $\shS_1^{q} \subsetdots \shS_t^q = \shS^{q}$ denote the strong Harder-Narasimhan filtration
of $\shS$.
Let $r_k = \rk (\shS^q_k/ \shS^q_{k-1})$,
$\bar{\mu}_k =   \mu(\shS_k^q/ \shS_{k-1}^q)/q$ and $\nu_k = - \bar{\mu}_k/\deg(Y)$,
$k=1 \komdots t$.
Write $\lceil q \nuk \rceil = q \nuk + \perk $ with the eventually periodic functions $\perk = \perk(q)$.
Let $ \nuu \leq  \nu_1$ and $  \nuo \gg \nu_t $ denote rational numbers and
set $  \lceil q \nuo \rceil =  q \nuo + \pero $.
Then for $q=p^{e} \gg 0$ we have $\sumceilnu {\shS} =$
\begin{eqnarray*}
\!& =& \frac{q^2}{2 \deg(Y)} 
\big( \mumu(\shS) + 2 \nuo \deg(\shS) \deg(Y) + \nuo^2 \rk(\shS) \deg(Y)^2 \big) \cr
\!& & + q (\nuo \rk(\shS) + \frac{\deg(\shS)}{\deg(Y)} ) (1-g - \frac{\deg(Y)}{2} )\cr
\!& & + q \pero ( \deg(\shS) + \nuo  \rk(\shS)\deg(Y)) \cr
\!& & + \rk(\shS)\pero  \big( (\pero -1) \frac{ \deg(Y)}{2} +1-g\big) 
- \sum_{k=1}^t r_k \perk \big( (\perk -1)\frac{ \deg(Y)}{2}+ 1-g\big)\cr
\!& &+ \sum_{k=1}^{t}
\big(\sum_{m = \lceil q \nu_k  \rceil }^{\lceil q \nu_k  \rceil  + \lceil\frac{\deg (\omega)}{\deg(Y)} \rceil}
h^1( (\shS_k/\shS_{k-1})^q(m)) \big)
\end{eqnarray*}
\end{theoremintro}
We emphasize that the right hand side is a simplification of the left hand side.
It is up to the last $h^1$-term a quadratic polynomial in $q$, where the linear term and the constant term have
eventually periodic coefficients determined by the eventually periodic functions $\pero$ and $\perk$.
In general we only know that the $h^1$-term is a bounded function of $q$
(Lemma \ref{bounded}). If however the groundfield
$K$ is the algebraic closure of a finite field, then this term is also eventually periodic
(Theorem \ref{periodic}).
This is due to the fact that the degree of the occurring
strongly semistable quotient sheaves
$(\shS_k/\shS_{k-1})^q(m)$ vary in a finite range, hence they
form a bounded family. Since they are defined over a finite field they behave eventually periodically in $q$.

In section \ref{exactsection} we look at a short exact sequence $0 \ra \shS \ra \shT \ra \shQ \ra 0$
and consider the alternating sum
$$ \sumceilnu \shS - \sumceilnu \shT +\sumceilnu \shQ \, $$
for suitable $\nuu$ and $\nuo$ (or summing over $\ZZ$).
We shall see that this equals
$$\frac{q^2}{2 \deg(Y)}(\mumu(\shS)- \mumu(\shT)+ \mumu(\shQ)) + O(q^0)$$
and that the $O(q^0)$-term is eventually periodic if everything is defined over a finite field (Theorem \ref{exactsequencetheorem}). This result underlines the significance
of the Hilbert-Kunz slope and shows that the expression
$(\mumu(\shS)-\mumu(\shT)+\mumu(\shQ))/2 \deg(Y)$
is an important invariant of a short exact sequences, which might be called its
Hilbert-Kunz multiplicity.

In the last section \ref{hilbertkunzsection} we apply this to the short exact syzygy sequence
$0 \ra \Syz(f_1 \komdots f_n)(0) \ra \bigoplus_{i=1}^n \O(-d_i) \ra \O(0) \ra 0$
given by ideal generators $f_i$ and obtain Theorem \ref{hilbertkunzfunctionperiodicintro}
(Theorem \ref{hilbertkunzfunctionperiodic}).
The Hilbert-Kunz multiplicity of this sequence equals the Hilbert-Kunz multiplicity of the ideal.

\section{The Hilbert-Kunz slope in positive characteristic}
\label{hilbertkunzslopesection}

Let $Y$ denote a smooth projective curve over an algebraically closed field $K$.
We recall briefly some notions about vector bundles, see  \cite{huybrechtslehn} for details.
The degree of a locally free sheaf $\shS$ on $Y$ of rank $r$
is defined by $\deg (\shS) = \deg \bigwedge^r (\shS)$.
The slope of $\shS$, written $\mu(\shS)$,
is defined by $\deg(\shS)/ r$.
The slope has the property that $\mu(\shS \otimes \shT)= \mu(\shS) + \mu(\shT)$.

A locally free sheaf $\shS$ is called semistable if
$\mu(\shT) \leq \mu(\shS)$ holds for every locally free subsheaf $\shT \subseteq \shS$.
Dualizing and tensoring with an invertible sheaf does not affect this property.

For every locally free sheaf $\shS$ on $Y$
there exists the so-called Harder-Nara\-sim\-han filtration
$\shS_1 \subsetdots \shS_t =\shS$,
where the $\shS_k$ are locally free subsheaves. This filtration is unique and has the property that
the quotients $\shS_{k}/ \shS_{k-1}$
are semistable and
$\mu (\shS_{k}/ \shS_{k-1}) > \mu(\shS_{k+1}/ \shS_{k})$ holds.

The number $\mu(\shS_1) = \mu_{\rm max}(\shS)$ is called the maximal slope of $\shS$,
and the minimal slope of $\shS$ is
$ \mu_{\rm min}(\shS) =\mu(\shS/ \shS_{t-1})$.
The existence of global sections can be tested with the maximal slope: if $\mu_{\max}(\shS)<0$,
then $H^0(Y, \shS)=0$.
Furthermore we have the relation $\mu_{\max} (\shS)= - \mu_{\min}(\shS^\dual)$,
where $\shS^\dual$ denotes the dual bundle.

With the help of the Harder-Narasimhan filtration we define in characteristic
zero the Hilbert-Kunz slope of a locally free sheaf by
$$ \mumu(\shS) = \sum_{k=1}^t r_k\mu_k^2 \, ,$$
where $r_k = \rk(\shS_k /\shS_{k-1})$ and $\mu_k= \mu(\shS_k /\shS_{k-1})$.
See \cite{brennerhilbertkunzcriterion} for the basic properties of this notion and its relation
to solid closure.

In positive characteristic we need the strong Harder-Narasimhan filtration.
We denote the pull-back of $\shS$  under the absolute Frobenius $F^{e}: Y \ra Y$
by $\shS^q$, $q = p^{e}$.
A locally free sheaf $\shS$ is called strongly semistable if $\shS^q$
is  semistable for every $q$.
Due to a theorem of A. Langer \cite[Theorem 2.7]{langersemistable}
there exists a Frobenius power
such that the quotients in the Harder-Narasimhan
filtration of the pull-back $\shS^q$ are all strongly semistable.
We call such a filtration the strong Harder-Narasimhan filtration and denote it by
$$ 0 \subset \shS_1^q \subsetdots \shS_t^q = \shS^q \, .$$
For $q' \geq q \gg 0$ the Harder-Narasimhan filtration of
$\shS^{q'}$ is
$$ \shS_1^{q'}= (\shS_1^q)^{q'/q} \subsetdots (\shS_t^{q'})^{q'/q} = (\shS_t^q)^{q'/q} \, .$$
This allows to define rational numbers
$ \bar{\mu}_k = \bar{\mu}_k (\shS) = \frac{\mu(\shS_k^q/ \shS_{k-1}^q)}{q}$ for $q \gg 0$.
The length $t$ of the strong Harder-Narasimhan filtration
as well as the ranks $r_k= \rk( \shS^q_k/ \shS^q_{k-1})$ are independent of $q \gg 0$.
For $q \gg 0$ we have $\mu (\shS_k^q/ \shS^q_{k-1}) =q \bar{\mu}_k$.

We now define the Hilbert-Kunz slope in positive characteristic.

\begin{definition}
\label{hilbertkunzslopedefinition}
Let $\shS$ denote a locally free sheaf on a smooth projective curve over an algebraically closed field of positive characteristic. Let $ \shS_1^q \subsetdots \shS_t^q = \shS^q$ denote
the strong Harder-Narasimhan filtration of $\shS$. Then the Hilbert-Kunz slope of $\shS$
is
$$ \mumu(\shS) = \sum_{k=1}^t r_k \bar{\mu}_k^2 \, .$$
\end{definition}

The Hilbert-Kunz slope is a rational number.
With this notion we may express the formula for the Hilbert-Kunz multiplicity
in the following way (\cite[Theorem 3.6]{brennerhilbertkunz}).

\begin{theorem}
Let $R$ denote a two-dimensional standard-graded normal domain
and let $I=(f_1 \komdots f_n)$ denote a homogeneous $R_+$-primary ideal
generated by homogeneous elements $f_i$ of degree $d_i, i=1 \komdots n$.
Then the Hilbert-Kunz multiplicity $e_{HK}(I)$ equals
$$e_{HK} (I) =\frac{ 1}{2 \deg(Y)}
(\mumu \big(\Syz(f_1 \komdots f_n)(0))-\deg(Y)^2 \sum_{i=1}^n d_i^2 \big) \, .$$ 
\end{theorem}

We shall see in the next section that the Hilbert-Kunz slope controls in general
the quadratic behavior (the $q^2$ term)
of the sections $\sumceilnu {\shS}$,
where $\nuo$ and $\nuu$ are rational numbers.

\section{Dimension of sections}
\label{dimensionsectionssection}

Lets fix the notation for this and the following sections.
Let $K$ denote an algebraically closed field of positive characteristic $p$.
Let $Y$ denote a smooth projective curve over $K$ of genus $g$
and canonical sheaf $\omega$. We fix a very ample invertible sheaf $\O(1)$ and denote by
$\deg(Y) = \deg(\O(1))$ the degree of the curve.
As usual we set $\shS(m)= \shS \otimes \O(m)$
for a coherent sheaf $\shS$ on $Y$. We shall only consider
locally free sheaves.

We want to describe the sum $\sumceilnu \shS$ for a locally free sheaf $\shS$,
where $\nuu$ and $\nuo$ are rational numbers.
It will become clear during this paper
that we cannot avoid rational boundaries. In fact we should be lucky that only rational
boundaries occur.
We may write $\lceil q \nuo \rceil = q \nuo + \pero (q) $, where
$0 \leq \pero (q) <1$. We often shall write $\pero$ instead of $\pero(q)$.
 If $\nuo = a/b$, then
$\lceil \frac{p^{e} a}{b} \rceil = \frac{p^{e} a}{b} + \frac{(- p^{e}a) \modu b}{b}$,
$0 \leq (- p^{e}a) \modu b < b$,
hence $\pero (q) =\frac{(- p^{e}a) \modu b}{b}$ is an eventually periodic function.
This is one source (the rounding source) of the periodic behavior of the Hilbert-Kunz function,
we treat the other  (the $H^1$ source) in section \ref{periodicitysection}.

\begin{lemma}
\label{sectionlemma}
Let $\shS$ denote a locally free sheaf on $Y$.
Let $\nuu < \nuo$ denote rational numbers.
Write $\lceil  q \nuu \rceil = q \nuu + \peru$
and $\lceil  q \nuo \rceil = q \nuo + \pero $.
Then
\begin{eqnarray*}
\sumceilnu {\shS} 
&=&\!\!\! q^2 \big(( \nuo- \nuu) \deg(\shS) 
+ ( \nuo^2 - \nuu^2) \frac{  \rk(\shS)\deg(Y)}{2}  \big) \cr
& &\!\!\!
+ q ( \nuo - \nuu) \rk(\shS) (1-g -\frac{\deg(Y)}{2} ) \cr
& &\!\!\! +q \big( ( \pero - \peru) \deg(\shS) + (\pero \nuo - \peru \nuu) \rk (\shS)\deg(Y) \big) \cr
& &\!\!\! + \rk (\shS) \big(  ( \pero(\pero\!-\!1) - \peru(\peru\!-\!1))  \frac{\deg(Y) }{2}
+ (\pero - \peru) (1\!- \!g) \big) \cr
& &\!\!\!+ \sumceilm h^1(\shS^q(m))
\end{eqnarray*}
\end{lemma}
\proof
Due to the formula of Riemann-Roch we have
\begin{eqnarray*}
\sumceilnu {\shS} & \!\!\! \!\! =\!\! \!\! & \! \! \! \sumceilm \! \! \deg(\shS^q(m))
\! + \! \! \sumceilm \! \rk(\shS) (1\!-\!g)
\!+\!\! \sumceilm \! h^1(\shS^q(m)) \cr
&=& \sumceilm  \big( q\deg(\shS) +m\rk(\shS) \deg(Y) \big)
\cr
& &+ ({\lceil q \nuo \rceil} -{\lceil q \nuu \rceil}) \rk(\shS) (1-g)
+ \sumceilm h^1(\shS^q(m)) \cr
&=& q ({\lceil q \nuo \rceil} -{\lceil q \nuu \rceil}) \deg(\shS)  \cr
& &
+\big( \lceil q \nuo \rceil ( \lceil q \nuo \rceil -1) - \lceil q \nuu \rceil (\lceil q \nuu \rceil - 1) \big)
\frac{ \rk(\shS)\deg(Y)}{2} \cr
& &+ ({\lceil q \nuo \rceil} -{\lceil q \nuu \rceil}) \rk(\shS) (1-g)
+ \sumceilm h^1(\shS^q(m)) \, .
\end{eqnarray*}
We now insert $\lceil q \nuu \rceil =  q \nuu + \peru (q)$ and $\lceil q \nuo \rceil =  q \nuo + \pero (q)$
in these summands. The first summand yields
$$q  (q(\nuo - \nuu) +\pero-\peru)\deg(\shS) =q^2 (\nuo - \nuu) \deg(\shS)
+q (\pero -\peru) \deg(\shS)\, .$$
The second summand yields
$$ \big( q^2 (\nuo^2-\nuu^2) -q(\nuo - \nuu)
+q(2\pero \nuo- 2 \peru \nuu)
+\pero(\pero-1)- \peru(\peru-1) \big) \frac{ \rk(\shS)\deg(Y)}{2} \, .$$
The third summand yields
$$ (q(\nuo-\nuu) +\pero-\peru) \rk(\shS) (1-g) \, .$$
We regroup the terms and get the $q^2$-term
$$ q^2 (\nuo - \nuu)\deg(\shS)  + q^2 (\nuo^2 - \nuu^2) \frac{ \rk(\shS)\deg(Y)}{2} \, ,$$
the constant $q$-term
$$-q (\nuo - \nuu) \frac{ \deg(Y) \rk(\shS)}{2} +q (\nuo - \nuu) \rk(\shS) (1-g) 
= q (\nuo - \nuu) \rk(\shS) (1-g- \frac{\deg(Y)}{2})
\, ,$$
the periodic $q$-term
$$ q (\pero - \peru)\deg(\shS) + q (2 \pero \nuo -2 \peru \nuu) \frac{ \rk(\shS)\deg(Y)}{2}$$
and the periodic $q^0$-term
$$  ( \pero(\pero -1) - \peru(\peru -1) )\frac{\rk(\shS)\deg(Y)}{2}
+ (\pero - \peru) \rk (\shS) (1-g) \, .$$
This is what we have written down.
\qed

\medskip
This Lemma is of course only useful if we can say something about the $h^1$-term.
We treat first the case of a strongly semistable sheaf $\shS$.

\begin{proposition}
\label{semistablesections}
Let $\shS$ denote a strongly semistable sheaf on $Y$.
Set $\nus = - \mu(\shS) /\deg(Y) = - \deg(\shS)/ \rk(\shS) \deg(Y)$.
Let $\nuu$ and $\nuo$ denote rational numbers such that
$ \nuu \leq \nus \ll \nuo$.
Let $ \lceil q \nus \rceil =  q \nus + \pers $ and
$\lceil q \nuo \rceil =  q \nuo + \pero $. Then
\begin{eqnarray*}
\sumceilnu {\shS} \!\!\!\!
& =& \!\!\!\! \frac{q^2}{2 \deg(Y)} 
\big( \mumu(\shS)\! +\! 2 \nuo \deg(\shS) \deg(Y)\! +\! \nuo^2 \rk(\shS) \deg(Y)^2 \big) \cr
& &\!\!\!\! + q  (\nuo \rk(\shS) + \frac{\deg(\shS)}{\deg(Y)})(1-g - \frac{\deg(Y)}{2} )  \cr
& & \!\!\!\! + q \pero  \big( \deg(\shS) + \nuo  \rk (\shS)\deg(Y) \big) \cr
& &\!\!\!\! + \rk(\shS) \big(   ( \pero(\pero -1) - \pers(\pers -1))\frac{\deg(Y) }{2}
+ (\pero - \pers) (1-g) \big) \cr
& & \!\!\!\!+ \sum_{m = \lceil q \nus  \rceil }^{\lceil q \nus  \rceil  + 
\lceil \frac{\deg (\omega)}{\deg(Y)} \rceil}
h^1(\shS^q(m))
\end{eqnarray*}
\end{proposition}
\proof
For $m < \lceil q \nus \rceil = \lceil - q \mu(\shS) /\deg(Y) \rceil$ we have
$m< -q \mu(\shS) / \deg (Y)$ and therefore
$$ \deg(\shS^q(m)) = q \deg(\shS) + m \rk(\shS)  \deg(Y) <  0\, .$$
Since $\shS^q$ is semistable we have $h^0(\shS^q(m))=0$ in this range.
Therefore we have
$$ \sum_{m = \lceil q \nuu   \rceil}^{\lceil q \nuo \rceil -1} h^0(\shS^q(m))
= \sum_{m = \lceil q \nus  \rceil}^{\lceil q \nuo \rceil -1} h^0(\shS^q(m)) \, .$$
We apply the formula from Lemma \ref{sectionlemma}
to the right hand side.
We insert $\nus =- \frac{\deg(\shS)}{\rk(\shS) \deg(Y)}$ and get
the quadratic term
\begin{eqnarray*}
& &  \deg(\shS) \big( \nuo + \frac{\deg(\shS)}{ \rk(\shS)\deg(Y)}\big)
+ \frac{\rk(\shS) \deg(Y)}{2} 
\big( \nuo^2 - \frac{\deg(\shS)^2}{\rk(\shS)^2 \deg(Y)^2} \big) \cr
&=& \nuo \deg(\shS) + \nuo^2 \frac{ \rk(\shS)\deg(Y)}{2}  
+ \frac{\deg(\shS)^2}{\rk(\shS) \deg(Y)}
- \frac{1}{2} \frac{ \deg(\shS)^2}{\rk(\shS) \deg(Y)} \cr
&=& \frac{1}{2 \deg(Y)}
\big( 2 \nuo \deg(\shS) \deg(Y)+ \nuo^2 \rk(\shS)\deg(Y)^2
+ \frac{ \deg(\shS)^2}{ \rk(\shS)} \big)\, .
\end{eqnarray*}
The linear-constant term yields
$$ \big( \nuo + \frac{\deg(\shS)}{ \rk(\shS) \deg(Y)}  \big)
 \rk(\shS) (1-g- \frac{\deg(Y) }{2}) \, .$$
The linear-periodic term is
\begin{eqnarray*}
&=&(\pero - \pers) \deg (\shS)
+ ( \pero \nus   + \frac{\pers \deg(\shS)}{\rk(\shS) \deg(Y)})  \rk(\shS) \deg(Y) \cr
&=& (\pero - \pers) \deg (\shS) + \pero \nus \rk(\shS)\deg(Y) + \pers  \deg(\shS)  \cr
&=& \pero \big(\deg(\shS) + \nus  \rk(\shS)\deg(Y) \big) \, .
\end{eqnarray*}
The forth term comes directly from the formula in Lemma \ref{sectionlemma}.
So now we have to look at
$\sum_{m = \lceil q \nus  \rceil}^{\lceil q \nuo \rceil -1} h^1(\shS^q(m))$.
For $m > \lceil q \nus \rceil+ \lceil \frac{\deg(\omega)}{\deg(Y)} \rceil$
we have
$$\mu( (\shS^q)^\dual (-m) \otimes \omega)
= -q \mu(\shS) -m \deg(Y) + \deg (\omega) < 0 \, ,$$
hence $h^1(\shS^q(m))=0$
by Serre duality and the semistability of $\shS^q$.
\qed

\section{Dimension of sections - general case}

We shall extend the results of the previous section to arbitrary locally
free sheaves using the strong Harder-Narasimhan filtration.
If $\shS_1^q \subsetdots \shS_t^q= \shS^q $
is the strong Harder-Narasimhan filtration of $\shS$,
we set $\bar{\mu}_k= \mu (\shS_k^q /\shS_{k-1}^q)/q$ and $\nu_k= -\bar{\mu}_k/ \deg(Y)$.

\begin{lemma}
\label{sectionsum}
Let $\shS$ denote a locally free sheaf on a smooth projective curve $Y$ over an algebraically closed field $K$ of positive characteristic $p$.
Let $\shS_k^{q} \subset \shS^{q}$ denote the strong Harder-Narasimhan filtration
of $\shS$. Then for numbers
$ \nuu \leq \nu_1$ and $  \nuo > \nu_t $ 
and for $q=p^{e} \gg 0$ we have
$$
\sum_{m= \lceil q \nuu  \rceil}^{\lceil q \nuo \rceil -1} h^0(\shS^q(m))
=  
\sum_{k=1}^t \big(  \sum_{m= \lceil q \nu_k \rceil}^{\lceil q \nuo \rceil -1}
h^0( (\shS_{k}/\shS_{k-1})^q (m))  \big) 
\, .
$$
\end{lemma}
\proof
We do induction on the strong Harder-Narasimhan filtration,
so assume that $q \gg 0$ such that the Harder-Narasimhan filtration
of $\shS^q$ is strong.
For $t=1$, that is in the strongly semistable case, the statement was proved in the beginning of the proof of
Proposition \ref{semistablesections}.
We use the short exact sequence
$0 \ra \shS_{t-1} \ra \shS \ra \shS/\shS_{t-1} \ra 0$
and we assume that the statement is true for $\shS_{t-1}$.
We write
$$ \sumceil {\nuu}{\nuo}{\shS} = \sumceil {\nuu} {\nut} {\shS} 
+ \sumceil  {\nut} {\nuo} {\shS} 
\, .$$
Look at the two summands on the right.
For the first summand we have
$m \leq \lceil q \nut \rceil -1 $,
hence $m < q \nut = -q \bar{\mu}_t /\deg(Y)$  and therefore
$$ \mu ( (\shS/ \shS_{t-1})^q (m)) = q \bar{\mu}_t + m \deg(Y) < 0 \, .$$
Hence $h^0((\shS/\shS_{t-1})^q (m) ) =0$ and therefore the short exact sequence yields
$\sumceil {\nuu} {\nut} {\shS} = \sumceil {\nuu} {\nut} {\shS_{t-1}}$.

Now look at the second summand. Since $ m \geq q \nut $
we have $ m > q\nu_{t-1}  +\frac{ \deg(\omega)}{\deg(Y)}$ for $q \gg 0$.
Therefore we have
$$ \mu_{\rm max} ( (\shS^q _{t-1})^\dual  (-m) \otimes \omega) 
= -q \bar{\mu}_{t-1} - m \deg(Y) + \deg(\omega) < 0$$
and so $H^1(\shS^q_{t-1}(m)) =0$ in this range.
Thus
$$\sumceil  {\nut} {\nuo} {\shS} = \sumceil  {\nut} {\nuo} {\shS_{t-1}}
+ \sumceil  {\nut} {\nuo} {(\shS/\shS_{t-1})} \, .$$
This gives the result.
\qed

\medskip
The following theorem describes the global sections of a locally free sheaf.

\begin{theorem}
\label{sectionformulatheorem}
Let $\shS$ denote a locally free sheaf on a smooth projective curve $Y$ over an algebraically closed field $K$ of positive characteristic $p$.
Let $\shS_1^{q} \subsetdots \shS_t^q = \shS^{q}$ denote the strong Harder-Narasimhan filtration
of $\shS$.
Let $r_k = \rk (\shS^q_k/ \shS^q_{k-1})$,
$\bar{\mu}_k =   \mu(\shS_k^q/ \shS_{k-1}^q)/q$ and $\nu_k = - \bar{\mu}_k/\deg(Y)$.
Write $\lceil q \nuk \rceil = q \nuk + \perk $ with the eventually periodic functions $\perk = \perk(q)$.
Let $ \nuu \leq  \nu_1$ and $  \nuo \gg \nu_t $ denote rational numbers and
set $  \lceil q \nuo \rceil =  q \nuo + \pero $.
Then for $q=p^{e} \gg 0$ we have $\sumceilnu {\shS} =$
\begin{eqnarray*}
\!& =& \frac{q^2}{2 \deg(Y)} 
\big( \mumu(\shS) + 2 \nuo \deg(\shS) \deg(Y) + \nuo^2 \rk(\shS) \deg(Y)^2 \big) \cr
\!& & + q (\nuo \rk(\shS) + \frac{\deg(\shS)}{\deg(Y)})(1-g - \frac{\deg(Y)}{2} ) \cr
\!& & +q \pero ( \deg(\shS) + \nuo \rk(\shS)\deg(Y)) \cr
\!& & + \rk(\shS)\pero  \big(  (\pero -1) \frac{\deg(Y)}{2} +1-g \big)
- \sum_{k=1}^t r_k \perk \big( (\perk -1) \frac{\deg(Y)}{2}+ 1-g\big)\cr
\!& &+ \sum_{k=1}^{t}
\big(\sum_{m = \lceil q \nu_k  \rceil }^{\lceil q \nu_k \rceil + \lceil \frac{\deg (\omega)}{\deg(Y)} \rceil}
h^1( (\shS_k/\shS_{k-1})^q(m)) \big)
\end{eqnarray*}
\end{theorem}
\proof
Let $q$ be big enough such that the statement of Lemma \ref{sectionsum} holds true,
so that we have
$$
\sum_{m= \lceil q \nuu  \rceil}^{\lceil q \nuo \rceil -1} h^0(\shS^q(m))
=  
\sum_{k=1}^t \big(  \sum_{m= \lceil q \nu_k \rceil}^{\lceil q \nuo \rceil -1}
h^0( (\shS_{k}/\shS_{k-1})^q (m))  \big) 
\, .$$
For $t=1$, that is in the strongly semistable case, the result is just Proposition
\ref{semistablesections} (with $\pers = \pi_1$).
In general we have to sum over the $t$ expressions coming from
Proposition \ref{semistablesections} 
for the strongly semistable quotient sheaves
$\shS_1^q, \shS_2^q/ \shS_1^q \komdots \shS^q/ \shS^q_{t-1}$.
The rank and the degree are additive on short exact sequences
and the Hilbert-Kunz slope is additive on the quotients in the strong Harder-Narasimhan filtration.
Therefore the summations of the first, the second and the third
terms from Proposition \ref{semistablesections} 
yield the first, the second and the third expression in the statement.
This is also true for the fifth term. For the forth term we just have to add
\begin{eqnarray*}
\!\!\!\!\!& & \!\!\!\!
\sum_{k=1}^t r_k \big(  ( \pero (\pero -1) - \perk (\perk -1)) \frac{\deg( Y) }{2}
+ (\pero - \perk )(1-g) \big) \cr
\!\!\!\!\!&=&\!\!\!\! \rk(\shS) \big(  \pero(\pero-1)\frac{\deg(Y)}{2} \!+\! \pero(1-g)\! \big)
\!-\! \sum_{k=1}^t r_k \big(  \perk (\perk -1) \frac{\deg(Y)}{2} + \perk(1-g)\! \big)
\end{eqnarray*}
\qed

\section{Boundedness and periodicity}
\label{periodicitysection}

We take now a closer look at the lower terms in the formula of Theorem \ref{sectionformulatheorem}.

\begin{lemma}
\label{bounded}
Let $\shS$ denote a strongly semistable locally free sheaf on the smooth projective curve
$Y$ over the algebraically closed field $K$ of positive characteristic $p$.
Set $\nus= - \frac{\deg(\shS)}{\rk(\shS) \deg(Y)}$ and let
$\nuo \geq \nus$ be a rational number.
Then $\sum_{m = \lceil q\nus \rceil}^{\lceil q \nuo \rceil} h^1(\shS^q(m)) = O(q^0)$.
\end{lemma}
\proof
It is enough to consider the sum running from
$\lceil q \nus \rceil $ 
to $\lceil q \nus \rceil + \lceil \frac{\deg (\omega)}{\deg(Y)} \rceil $,
since above this we have $H^1(Y, \shS^q(m))=0$.
Hence we are concerned with $\lceil \frac{\deg (\omega)}{\deg(Y)} \rceil+1 $ summands.
From $ \frac{-q \deg(\omega)}{ \rk(\shS) \deg(Y)} \leq \lceil q \nus \rceil
\leq m \leq \lceil q \nus \rceil + \lceil \frac{\deg(\omega)}{\deg(Y)} \rceil \leq
q \nus + \frac{\deg(\omega)}{\deg(Y)} +2 $
we get
$0 \leq \deg(\shS^q(m)) \leq \rk(\shS) (\deg(\omega)+ 2 \deg(Y))$.
Hence the degrees of these sheaves vary in a finite range.

It follows now from fundamental boundedness results for semistable sheaves
that there exists an upper bound for the dimension of global sections of these sheaves
(\cite[Corollary 1.7.7]{huybrechtslehn}). This  is then also true for $h^1( \shS^q(m))$.
\qed

\begin{remark}
In the construction via quot-schemes of the moduli space of semistable bundles of given degree
one shows that there exists a coherent sheaf $\shF$ such that every $\shS$ admits a surjection
$\shF \ra \shS \ra 0$. From this it follows again that
$h^1(\shS)$ is bounded.
For invertible sheaves of fixed degree 
one may use the theorem of Clifford, see \cite[Theorem IV.5.4]{haralg}, to show
that there exists a common bound for the dimension of their global sections.
Without the condition semistable this conclusion does not hold, as
the example $\shL \oplus \shL^{-1}$ shows.
\end{remark}

\begin{theorem}
\label{periodic}
Let $K$ denote an algebraically closed field of positive characteristic $p$,
let $Y$ denote a smooth projective curve over $K$.
Let $\shS$ denote a locally free sheaf on $Y$.
Let $ \nuu \leq \nu_1$ and $  \nuo \gg \nu_t $ denote rational numbers.
Then we have
$$  \sumceilnu {\shS}    = \alpha  q^2  + \beta(q) q + \gamma(q)  \, ,$$
where $\alpha = \frac{\mumu(\shS) +2\nuo \deg(\shS) + \nuo^2 \rk(\shS) \deg(Y)^2}{2 \deg(Y)}$
is a rational number,
$\beta(q)$ is an eventually periodic function and $\gamma(q)$ is a
bounded function {\rm(}both with rational values{\rm)}.

If moreover $K$ is the algebraic closure of a finite field,
then $\gamma(q)$ is also an eventually periodic function.
\end{theorem}
\proof
The first statement follows directly from Theorem \ref{sectionformulatheorem}
and Lemma \ref{bounded}.
For the second statement  we only have to show that the expression
$$\sum_{k=1}^{t}
\big(\sum_{m = \lceil q \nu_k  \rceil }^{\lceil q \nu_k  \rceil
+ \lceil \frac{\deg (\omega)}{\deg(Y)} \rceil }
h^1(Y, (\shS_k/\shS_{k-1} )^q(m)) \big) $$
is an eventually periodic function in $q$. Thus we may assume that $\shS$ is strongly semistable.

Set $\nus =- \frac{\deg(\shS)}{\rk(\shS) \deg(Y)}$.
We consider the starting term of the summation, $m(q) = \lceil q \nus \rceil$.
Write $m(q) = q \nus + \pert(q)$ with the eventually periodic function $\pert(q)$.
Let $\tilde{q}$ be the length of the periodicity.
We have
$$\deg ( \shS^q(m(q)))= q \deg (\shS) + m(q) \rk(\shS) \deg(Y)
= \pert (q) \rk(\shS) \deg(Y) \, ,$$
so that the degree of these sheaves behaves also eventually periodical with
the same periodicity $\tilde{q}$.

We consider now a subset of type $M= \{q_0 \tilde{q}^\ell:\, \ell \in \NN \}$.
In particular $\pert (q)$ is constant on this set $M$.
For $q \in M$ we have
\begin{eqnarray*}
\shS^{q \tilde{q}}(m(q \tilde{q}))
&=& \shS^{q \tilde{q}} \otimes \O(q \tilde{q} \nus + \pert(q \tilde{q})) \cr
&=& (\shS^q)^{\tilde{q}} \otimes \O(q \tilde{q} \nus + \pert(q) +( \tilde{q} -1) \pert(q))
\otimes \O( (- \tilde{q}+1) \pert(q)) \cr
&=& (\shS^q)^{\tilde{q}}( q \tilde{q} \nus + \tilde{q} \pert(q)) \otimes \O( (- \tilde{q}+1) \pert (q)) \cr
&=& (\shS^q(q\nus + \pert(q)))^{\tilde{q}} \otimes \O( (- \tilde{q}+1) \pert (q))
\end{eqnarray*}
This means that the $\tilde{q}$ successor is build from its predecessor by pulling it back
and tensor the result with a fixed invertible sheaf. In particular this recursion rule is independent
of $q$.
Now the curve and the locally free sheaf $\shS$ are defined over a finite subfield $L \subset K$.
Then all the $ \shS^q (m(q))$, $q \in M$, have the same degree and are defined over $L$.
The family of semistable sheaves with fixed degree is bounded and therefore
there exist only finitely many such sheaves defined over $L$.
Hence the recursion above shows that the sequence of sheaves
$\shS^q (m(q))$, $q \in M$, is eventually periodic.

From this it follows immediately that also the other sheaves
$\shS^q(m(q) +s ) = \shS^q(m(q)) \otimes \O(s)$ for fixed
$s$, $ 0 \leq s \leq \lceil \frac{\deg (\omega)}{\deg(Y)} \rceil$,
occur periodically in $q \in N$.
This is then also true for their $h^1$-term.
\qed

\begin{remark}
A similar argument was used by Lange and Stuhler in \cite{langestuhler}
to show that the Frobenius pull-backs of
a strongly semistable sheaf of degree $0$ on a curve over a finite field
behave eventually periodically.
\end{remark}

\section{Short exact sequences}
\label{exactsection}

In this section we look at a short exact sequence
$0 \ra \shS \ra \shT \ra \shQ \ra 0$ of locally free sheaves
on  a smooth projective curve $Y$ over an algebraically closed field $K$
of positive characteristic $p$. We want to compute the alternating sum
$$\sum_{m \in \ZZ} \big(h^0(\shS^q(m)) - h^0(\shT^q(m))+ h^0(\shQ^q(m)) \big)$$
in dependence of $q=p^{e}$. 
We will see in the next section that the computation of the Hilbert-Kunz function of an ideal
in a two-dimensional normal standard-graded domain is a special case of this consideration.
The sum is for every $q$ finite, since for $m \ll 0$ all terms are $0$ and for $m \gg 0$
we have $H^1(Y, \shS^q(m))=0$ and the sum is $0$.
The sum is the dimension of the cokernel
$ \sum_{m \in \ZZ} \dim ( \Gamma(Y, \shQ^q(m)) / \im (\Gamma(Y, \shT^q(m) ))$
and equals also
$$\sum_{m \in \ZZ} \big(h^1(\shS^q(m)) - h^1(\shT^q(m))+ h^1(\shQ^q(m)) \big) \, .$$
The following theorem shows that this alternating sum is a quadratic polynomial
with the alternating sum of the Hilbert-Kunz slopes as leading coefficient,
with vanishing linear term and with a bounded constant term, which is periodic in the finite case.

\begin{theorem}
\label{exactsequencetheorem}
Let $Y$ denote a smooth projective curve over an algebraically closed field
of positive characteristic $p$.
Let $0 \ra \shS \ra \shT \ra \shQ \ra 0$
denote a short exact sequence of locally free sheaves on $Y$.
Then the following hold.
\numiii
\begin{enumerate}
\item
The alternating sum of the global sections is
\begin{eqnarray*}
& &\sum_{m \in \ZZ} \big(h^0(\shS^q(m)) - h^0(\shT^q(m))+ h^0(\shQ^q(m)) \big) \cr
&=&\frac{q^2}{2 \deg(Y)} (\mumu(\shS) -\mumu(\shT)+ \mumu(\shQ)) + O(q^0) \, .
\end{eqnarray*}
\item
Let $\shS_k \subseteq \shS$, $\shT_i \subseteq \shT$ and
$\shQ_j \subseteq \shQ$ denote the strong Harder-Narasim\-han filtrations of these sheaves
with disjoint index sets.
Let $r_k$, $\bar{\mu}_k$, $\nu_k$ and $\pi_k$ {\rm(}$r_i$, $r_j$, $\nu_i$, $\nu_j$ etc.  respectively{\rm)}
denote the ranks, slopes and the eventually
periodic functions corresponding to the strongly semistable quotients.
Then the $O(q^0)$-term equals
\begin{eqnarray*}
\! \!\!\!\!\! & &\! \!- \!\sum_{k} r_k \pi_k \big(( \pi_k-1)\frac{ \deg(Y)}{2}+1-g \big)
\!+\!\sum_{i} r_i \pi_i \big( ( \pi_i-1)\frac{ \deg(Y)}{2}+1-g \big) \cr
\!\! \!\!\! \! & &\!\!-\! \sum_{j} r_j \pi_j \big( ( \pi_j-1)\frac{ \deg(Y)}{2}+1-g\big) 
\! +\!\sum_k \big( \sum_{m = \lceil q \nu_k \rceil}^{\lceil q \nu_k \rceil + 
\lceil \frac{\deg(\omega)}{\deg(Y)} \rceil}
h^1(  (\shS_k^q/\shS_{k-1}^q) (m))\big) \cr
\!\!\!\!\!\! & &\!\!
-\! \sum_i \big(\! \sum_{m = \lceil q \nu_i \rceil}^{\lceil q \nu_i \rceil + \lceil  \frac{\deg(\omega)}{\deg(Y)} \rceil}
  h^1(  (\shT_i^q/\shT_{i-1}^q) (m))\big) 
\!+\!\sum_j \big(\!
\sum_{m = \lceil q \nu_j \rceil}^{\lceil q \nu_j \rceil +\lceil \frac{\deg(\omega)}{\deg(Y)} \rceil}
h^1(  (\shQ_j^q/\shQ_{j-1}^q) (m))\big) .
\end{eqnarray*}

\item
If $K$ is the algebraic closure of a finite field, then
the $O(q^0)$-term is eventually periodic.
\end{enumerate}
\end{theorem}

\proof
Let $\nuu$ and $\nuo$ be rational numbers such that
$\nuu \leq  \nus_i, \nus_j, \nus_k \ll \nuo$ for all $k,i,j$.
Then we may look at the finite alternating sum running from
$\lceil q \nuu \rceil $ to $\lceil q \nuo \rceil $.
We only have to add the expressions in Theorem \ref{sectionformulatheorem} for $\shS$, $\shT$ and $\shQ$.
Since the rank and the degree are additive on short exact sequences most terms vanish.
The remaining terms are
$\frac{q^2}{2\deg(Y)}( \mumu(\shS)- \mumu(\shT)+ \mumu(\shQ))$
and the terms written down in (ii).
The boundedness stated in (ii) follows from Lemma \ref{bounded}
and the periodicity statement in
(iii) follows from Theorem \ref{periodic}.
\qed

\section{The Hilbert-Kunz function of an ideal}
\label{hilbertkunzsection}

We come back to the Hilbert-Kunz function of an ideal.
Let $R$ denote a normal standard-graded
domain over an algebraically closed field of positive characteristic $p$
and let $f_1 \komdots f_n$ denote homogeneous generators
of an $R_+$-primary ideal of degrees $d_1 \komdots d_n$.
These data give rise to the short exact sequence of locally free sheaves on
$Y= \Proj R$,
$$ 0 \lra \Syz(f_1 \komdots f_n)(m) \lra
\bigoplus_{i=1}^n  \O(m-d_i) \stacklra {f_1 \komdots f_n}  \O(m) \lra 0 \, .$$
The pull-back of this short exact sequence
under the $e$-th absolute Frobenius morphism $F^{e}: Y \ra Y$ yields
$$ 0 \lra (F^{e}(\Syz(f_1 \komdots f_n))) (m) \lra
\bigoplus_{i=1}^n  \O(m-qd_i) \stacklra {f_1^q \komdots f_n^q}  \O(m) \lra 0 \, .$$
Since $R$ is normal the global sections $\Gamma(Y, -)$ of this sequence
yield
$$ 0 \lra \Gamma(Y, \Syz (f_1^q \komdots f_n^q)(m))
\lra \bigoplus_{i=1}^n R_{ m-qd_i}
\stackrel{f_1^q \komdots f_n^q} {\lra} R_m  \lra \ldots $$
and the cokernel of the last mapping is just $(R /I^{[q]}) _m$.
Hence we have 
$$ \length  ((R/I^{[q]}) _m)
= h^0 (\O(m)) - \sum_{i=1}^n  h^0(\O(m-qd_i ))+ h^0(\Syz (f_1^q \komdots f_n^q)(m) )  \, . $$
Therefore we may express the Hilbert-Kunz function
$$ \hkf (q) = \sum_{m=0}^\infty \length (R/I^{[q]})$$
as the alternating sum of this short exact sequence.
Hence we apply the results of the previous section to this situation
and deduce the following theorem.

\begin{theorem}
\label{hilbertkunzfunctionperiodic}
Let $K$ denote an algebraically closed field of positive characteristic $p$.
Let $R$ denote a normal two-dimensional standard-graded
$K$-domain and let $I$ denote a homogeneous
$R_+$-primary ideal.
Then the Hilbert-Kunz function of $I$ has the form
$$ \hkf (q) =e_{HK}(I) q^{2} + \gamma(q) \, ,$$
where $e_{HK}(I)$ is a rational number and $\gamma(q)$ is a bounded function.

Moreover, if $K$ is the algebraic closure of a finite field,
then $\gamma(q)$ is an eventually periodic function.
\end{theorem}
\proof
This follows at once from Theorem \ref{exactsequencetheorem}
applied to the short exact syzygy sequence.
\qed

\begin{remark}
The rationality of the leading coefficient was proved in \cite[Theorem 3.6]{brennerhilbertkunz}.
That the Hilbert-Kunz function has the form $e_{HK}(I) q^{2} +  \beta q +\gamma(q)$
with real numbers $e_{HK}(I)$ and $\beta$ and a bounded function $\gamma(q)$
was proved in \cite{hunekemcdermottmonsky}. In that paper Monsky also announces that
$\beta =0$. The periodicity was only known for the maximal ideal in the elliptic case (\cite{buchweitzchenhilbertkunz},
\cite{monskyellipticchartwo}, \cite{fakhruddintrivedi}), where the finiteness condition is not necessary. 
The periodicity in general for an arbitrary field is completely open.
\end{remark}

\bibliographystyle{plain}

\bibliography{bibliothek}

\end{document}